\documentclass[a4paper,12pt,fleqn,leqno]{article}

\usepackage[utf8]{inputenc}
\usepackage[greek,english]{babel}
\usepackage{amssymb,amsmath,amsthm}
\usepackage{tikz}
\usepackage[colorlinks=true,allcolors=blue]{hyperref}

\newcommand {\comment}[1] {}
\newcommand {\df}[1] {{\bfseries #1}}
\newcommand {\qt}[1] {``{#1}''}
\renewcommand {\`}{\grave{\ }}
\newcommand {\area} {{\textstyle\int}}

\newcommand {\cC} {{\cal C}}
\newcommand {\cD} {{\cal D}}
\newcommand {\cE} {{\cal E}}
\newcommand {\cK} {{\cal K}}
\newcommand {\cL} {{\cal L}}
\newcommand {\cM} {{\cal M}}
\newcommand {\cX} {{\cal X}}

\newcommand {\stride} {\text{\scalebox{-0.4}[1]{$\angle$}}}
\newcommand {\enum}[2] {\{{#1},\dots,{#2}\}}
\newcommand {\eps} {\varepsilon}

\newcommand {\id}{{\mathrm{id}}}
\newcommand {\mul} {{\:\!\circ\!\!\!\!\;\cdot\,}}
\newcommand {\N} {{\mathbb N}}

\newcommand {\R} {{\mathbb R}}
\newcommand {\rest}[1] {|_{#1}}

\newcommand {\set}[2] {\{#1\sothat#2\}}
\newcommand {\sm}{\setminus}

\newcommand {\sothat} {\,:\,}
\newcommand {\sseq} {\subseteq}
\newcommand {\ti}[1] {\tilde{#1}}
\newcommand {\thet} {\vartheta}
\newcommand {\thru}[3] {{#1}_{#2},\dots,{#1}_{#3}}

\newcounter{substate}
\renewcommand{\thesubstate}{(\alph{substate})}
\newenvironment{substate}
{
  \begin{list}{\bf\thesubstate}
    {\usecounter{substate}
     \itemindent0em
     \settowidth\labelwidth{\bf(g)} \labelsep0.5em  
     \leftmargin\labelwidth \addtolength\leftmargin\labelsep
     \topsep0.5ex
     \itemsep0ex}
}
{\end{list}}

\newenvironment {subproof}
{
  \begin{list}{\bf\thesubstate}
    {\usecounter{substate} 
     \leftmargin0em 
     \settowidth\labelwidth{\bf(a)} \labelsep0.5em  
     \itemindent\labelwidth \addtolength\itemindent\labelsep
     \topsep0.5ex
     \itemsep0ex}
}
{\end{list}}

\newenvironment{mycases}
{
  \begin{list}{$\bullet$}
    {\itemindent0em
     \settowidth\labelwidth{$\bullet$}\labelsep0.5em  
     \leftmargin\labelwidth\addtolength\leftmargin\labelsep
     \topsep0.5ex
     \itemsep0ex}
}
{\end{list}}

\swapnumbers
\newtheorem{thm}{Theorem}[section]
\newtheorem{rmk}[thm]{Remark}
\newtheorem{lem}[thm]{Lemma}
\newtheorem{stridelem}[thm]{Stride Lemma}
\newtheorem{hammocklem}[thm]{Hammock Lemma}
\newtheorem{prp}[thm]{Proposition}
\newtheorem{cor}[thm]{Corollary}

\begin{document}
\title{Solution to the iterative differential
       equation $-\gamma g' = g^{-1}$}
\author{Roland Miyamoto}

\maketitle

\begin{abstract}
\noindent
Using a Picard-like operator $T$,
we prove that the iterative differential equation
$-\gamma g' = g^{-1}$ with parameter $\gamma>0$
has a solution $g=h\colon[0,1]\to[0,1]$
for only one value $\gamma=\kappa\approx0.278877$,
and that this solution $h$ is unique.
As an even stronger result,
we exhibit $h$ as the global limit
of the operator $T$.
\end{abstract}



\section{Introduction}

The study of Levine's sequence~\cite{Le,MaPoSl,Sl}
naturally leads to a certain differentiable function
$g\colon[0,1]\to[0,1]$
satisfying the following property:
When we rotate (the graph of) $g$ clockwise
by $90^\circ$ about the origin and
subsequently stretch it vertically
by a suitable positive factor,
then we obtain its derivative~$g'$.
Formally speaking, $g$ satisfies
the iterative differential equation (IDE)
\begin{equation}\label{eqn-ide}
  g\colon[0,1]\to[0,1],\quad
  -\gamma g' = g^{-1}\quad\text{for some}\quad\gamma>0
\end{equation}
where $g^{-1}$ denotes the compositional inverse of~$g$.
We call such a function a \df{unit stribola}
(from Greek \textgreek{στρίβω} = turn, twist).
Every unit stribola, that is, any solution $g$ to~\eqref{eqn-ide}
will obviously be continuously differentiable
and strictly decreasing and satisfy the identities
\[
  \textstyle
  g(0)=1,\quad g(1)=0,\quad
  g'(0)=-\frac1\gamma,\quad
  \int_0^1g=\gamma.
\]
IDEs similar to~\eqref{eqn-ide} have been studied by
Eder~\cite{Ed}, Fe\v{c}kan~\cite{Fe}, Buic\u{a}~\cite{Bu},
Egri and Rus~\cite{EgRu} and Berinde~\cite{Be},
but the techniques employed there
appear not to suit our situation.

\begin{figure}\stepcounter{figure}\label{fig-h}
\begin{center}
\begin{tikzpicture}
  \draw (0, 0) node[inner sep=0]
  {\includegraphics[width=0.36\textwidth]{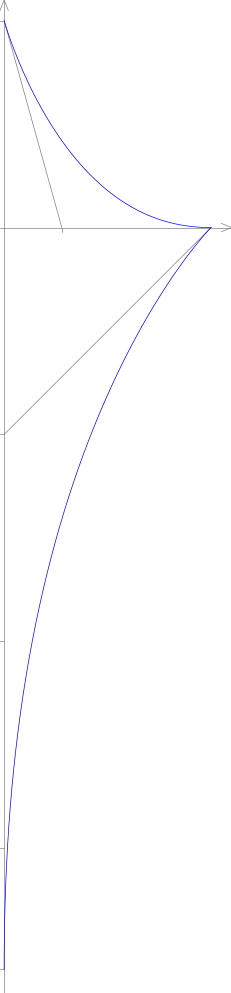}};
  \draw (-0.2,  8.0 ) node{\color{blue}$h$};
  \draw (-0.3,  1.7 ) node{\color{blue}$h'$};
  \draw (-1.32, 6.38) node{$\kappa$};
  \draw (-3.1, 11.88) node{$1$};
  \draw (-3.1,  6.7 ) node{$0$};
  \draw (-3.25, 1.56) node{$-1$};
  \draw (-3.3, -3.65) node{$-2$};
  \draw (-3.3, -8.84) node{$-3$};
  \draw (-3.3,-11.85) node{$-\frac1\kappa$};
  \draw (-2.2,  7.7 ) node{$\frac\kappa2$};
  \draw (-0.9,  7.4 ) node{$\frac\kappa2$};
  \draw (-1.1,  4.5 ) node{$\frac12$};
  \draw (-2.0,  0.3 ) node{$\frac12$};
  \draw ( 2.14, 0.55) node{Figure~\thefigure:};
  \draw ( 5.0,  0.0 ) node{The unit stribola {\color{blue}$h$}
                           and its derivative {\color{blue}$h'$}.};
  \draw ( 4.72,-0.55) node{Their tangents at $0$ resp.~$1$ bisect the};
  \draw ( 4.35,-1.02) node{areas they enclose with the axes.};
  \draw ( 4.60,-1.65) node{The derivative $h'$ arises from $h$ by a};
  \draw ( 4.90,-2.15) node{$90^\circ$ clockwise rotation about the origin};
  \draw ( 4.65,-2.63) node{and a subsequent vertical stretching};
  \draw ( 4.54,-3.10) node{with $\frac1\kappa$ where
                           $\kappa=\int_0^1h\approx0.278877$.};
\end{tikzpicture}
\end{center}
\end{figure}

In~\cite{MiSa}, we have constructed
a unit stribola $h$ by an iterative process.
At each step of this process, we perform the following
operation, denoted $T$:
Given any decreasing function $f\colon[0,1]\to[0,1]$
with $f(0)=1$ and area $\alpha:=\int_0^1f \neq 0$,
we rotate it clockwise by $90^\circ$ about the origin,
then stretch it vertically by $\frac1\alpha$,
then integrate, to obtain $Tf\colon[0,1]\to[0,1]$.
Starting from the line segment $h_1=1-\id_{[0,1]}$,
the sequence of iterates $h_1,\,h_2:=Th_1,\,h_3:=Th_2,\dotsc$
is shown to converge to a unit stribola~$h$.
Figure~\ref{fig-h} on the next page illustrates the functions
$h$ and $h'$ and some of their geometric properties.

\pagestyle{empty}
\newpage
\pagestyle{plain}

For the sake of better clarity,
we include the existence proof from~\cite{MiSa}
here in a simplified form.
A key ingredient to proving the convergence $h_n\to h$
is the observation
that the $h_n$, when stretched horizontally and vertically by arbitrary positive factors,
always \qt{cross} each other at most twice.
This feature (addressed by the concept of \qt{domination})
automatically extends to their limit.
To be explicit,
the unit stribola $h$ is \qt{dominated} by each iterate $h_n$.
From this, we can easily prove $h$ to be the only unit stribola.
Finally, we will establish global convergence, that is,
$\lim_{n\to\infty}T^nf = h$
for any decreasing function $f\colon[0,1]\to[0,1]$
with $f(0)=1$ and $\int_0^1f > 0$.
We attain this strong result by again and more heavily exploiting
the domination structure between the iterates~$h_n$ and their limit~$h$.

\section{The operator $T$}

For $0\leq a\leq b\leq1$
and any (Lebesgue) measurable function
$f\colon [0,1]\to[0,\infty)$,
we abbreviate $\int_a^bf := \int_a^bf(x)dx$
and $\area f :=\int_0^1f$.
We will also conveniently write
$\id:=\id_{[0,1]}$
for the identity function on~$[0,1]$.
Our investigations will involve the spaces
\begin{align*}
  \cM &:= \set{f\colon [0,1]\to[0,\infty)}
            {f\text{ measurable, }\area f>0},\\
  \cE &:= \set{f\in\cM}
            {f\text{ decreasing, }f(0)=1},\\
  \cC &:= \set{f\in\cE}{f\text{ continuous},\,f(1)=0},
  &\breve\cC &:= \set{f\in\cC}{f\text{ convex}},\\
  \cD &:= \set{f\in\cC}{f\text{ strictly decreasing}},
  &\breve\cD &:= \cD\cap\breve\cC,\\
  \cD' &:= \set{f\in\cD}{f\text{
          continuously differentiable on
          }(0,1]},
  &\breve\cD' &:= \cD'\cap\breve\cC,\\
  \cD\` &:=\set{f\in\cD'}
              {f'(1)=0,\;
                \lim_{x\to0}f'(x)\in(-\infty,0]
                \text{ exists}},
  &\breve\cD\` &:= \cD\`\cap\breve\cC
\end{align*}
of functions on $[0,1]$.
On $\cE$ we consider the \df{$\sup$-metric} $d_\infty$
and the \df{$1$-pseudometric} $d_1$ defined by
\[
  d_\infty(f,g) := \sup_{x\in[0,1]}|f(x)-g(x)|
  \quad\text{and}\quad
  d_1(f,g) := \area|f-g|
  \quad\text{for } f,g\in\cE.
\]
Note that $d_1$ is not a metric on $\cE$, but on $\cC$
because $d_1(f,g)=0 \iff f=g$ almost everywhere.
With any given $g\in\breve\cC$
we associate its \df{stride}
\begin{equation}\label{eqn-stride}
  \stride g := \sup\set{\alpha\geq0}{\alpha-\id\leq\alpha g}
  \in[0,1]
\end{equation}
and note that the two slopes
\begin{equation}\label{eqn-slopes}
  g'(0):= -\tfrac1{\stride g}\in[-\infty,-1]
  \quad\text{and}\quad
  g'(1) = \inf\set{\alpha\leq0}{\alpha(\id-1)\leq g}
          \in [-1,0]
\end{equation}
are well-defined, where we intentionally allow $g'(0)$
to assume the value $-\infty$.

Given $g\in\cE$,
we will use its \df{pseudo-inverse}
$g^*\in\cE$ defined by
\[
  g^*(y) := \sup g^{-1}[y,1]
          = \sup\set{x\in[0,1]}{g(x)\geq y}
  \quad\text{for } y\in[0,1].
\]
According to~\ref{rmk-inv}(c) below,
$g^*$ equals the compositional inverse $g^{-1}$ if $g\in\cD$.
Thus we may, and shall, consistently write $g^*$
in all cases from now on.
In~\cite{MiSa}, the following properties of the pseudo-inverse
are established.
\begin{rmk}\label{rmk-inv}
  For $f,g\in\cE$, the following statements hold.
  \begin{substate}
  \item $f\leq g \implies f^*\leq g^*$.
  \item If $g\in\cC$, then $g^*$ is strictly decreasing.
  \item If $g\in\cD$, then $g^*=g^{-1}\in\cD$
    is the inverse function of~$g$.
  \item $g\in\breve\cD \iff g^*\in\breve\cD$.
  \item $g\in\breve\cD\` \implies g^*\in\breve\cD'$.
  \item $\area g=\area g^*$.
  \item $d_1(f,g) = \area|f-g| = \area|f^*-g^*| = d_1(f^*,g^*)$.
  \item $\int_{g(t)}^1g^* = \int_0^tg - tg(t)$
    for all $t\in[0,1]$.
  \end{substate}
\end{rmk}

For given $g\in\cM$ and $f\in\breve\cD\`$,
we define the continuous functions
$Ig$ and $Df$ by setting
\[
  (Ig)(x) := \frac{\int_x^1g}{\area g}
  \quad\text{and}\quad
  (Df)(x):=\frac{f'(x)}{f'(0)}
  \quad\text{for } x\in[0,1]
\]
and formally introduce the operator $T\colon\cE\to\cC$
described in the introduction by setting
\[
  Tf := If^* \text{ for } f\in\cE,
  \quad\text{that is,}\quad
  (Tf)(x) = \frac{\int_x^1f^*}{\area f}
  \text{ for } x\in[0,1]
\]
by~\ref{rmk-inv}(f), as well as its iterations
$T^0=\id_\cE$ and $T^n:=T\circ T^{n-1}$ for $n\in\N$.

\begin{prp}\label{prp-DIT}
  For $g\in\cC$,
  the following statements hold.
	\begin{substate}
	\item $Tg\in\breve\cD\sm\{1-\id\}$.
	\item $1-\frac\id{{\int}g} \leq Tg \leq 1-\id$.
	\item $(Tg)(g(t)) \cdot\area g = \int_0^tg-tg(t)$
    for all $t\in[0,1]$.
  \item If $g\rest{[0,1)}>0$, then
    $Ig\in\breve\cD\`$
    with $-\frac1{(Ig)'(0)} = \stride Ig = \area g$
    and $DIg=g$.
  \item If $g\in\breve\cD\`$, then $Dg\in\cC$
    with $(Dg)\rest{[0,1)}>0$,
    $\area Dg = \stride g$ and $IDg=g$.
	\item $\cC^>:=\set{f\in\cC}{f\rest{[0,1)}>0}\sseq\cD$,
	  and $I\rest{\cC^>}\colon\cC^>\to\breve\cD\`$ is bijective
	  with inverse $D$.
	\end{substate}
\end{prp}
\begin{proof}
  \begin{subproof}
  \item Let $0\leq a<b\leq1$.
	  Then $g^*(a)>g^*(b)$ by~\ref{rmk-inv}(b), hence
	  \[
	    \area g \cdot \bigl[(Tg)(a)-(Tg)(b)\bigr]
	    = \int_a^bg^* > (b-a)g^*(b) \geq 0,
	  \]
	  showing that $Tg$ is strictly decreasing.
	  Similarly, for $a<x<b$, we obtain
	  \begin{equation}\label{eqn-convex}
	  	  \area g\cdot\frac{(Tg)(a)-(Tg)(x)}{x-a}
	      = \frac{\int_a^xg^*}{x-a} > g^*(x)
	      > \frac{\int_x^bg^*}{b-x}
	      = {\area g}\cdot\frac{(Tg)(x)-(Tg)(b)}{b-x},
	  \end{equation}
	  hence $Tg$ is convex.
	  Evaluating~\eqref{eqn-convex} for $0=a<x<b=1$
	  yields $(Tg)(x)<1-x$,
	  which completes the proof of the assertion.
	\item We have $\gamma:=\area g = \area g^*$
	  according to~\ref{rmk-inv}(f)
	  and $g^*\leq1$, hence
	  $\gamma(Tg)(x) = \int_x^1g^* = \gamma-\int_0^xg^*
	   \geq \gamma-x$ for all $x\in[0,1]$,
	  implying the left inequality,
	  while the right one follows from~(a).
	\item follows from~\ref{rmk-inv}(f) and (h).
  \item By the assumptions,
    $Ig\in\cD$,
    $(Ig)' = -\frac g{{\int}g}$
    is continuous and increasing,
    $(Ig)'(1) = 0$,
    and we have
    $-\frac1{\stride Ig} = (Ig)'(0) = -\frac1{{\int}g}$
    using~\eqref{eqn-slopes},
    so that $Ig\in\breve\cD\`$
    and $DIg=\frac{(Ig)'}{(Ig)'(0)}=g$.
  \item From $g\in\breve\cD\`$ we conclude that
    $g'(0)<-1$ and
    $Dg=\frac{g'}{g'(0)}:[0,1]\to[0,1]$ is continuous
    and decreasing with $(Dg)\rest{[0,1)}>0$,
    $(Dg)(0)=1$ and $(Dg)(1)=0$, hence $Dg\in\cC$.
    Moreover,
    $\area Dg = \frac{g(1)-g(0)}{g'(0)} = \stride g$
    by~\eqref{eqn-slopes}
    and $(\area Dg)(IDg)(x) = \frac{\int_x^1g'}{g'(0)}
     = \stride g \cdot g(x)$ for all $x\in[0,1]$,
    hence $IDg=g$.
  \item follows from~(c) and~(d).
    \qedhere
  \end{subproof}
\end{proof}

We now explicitly state the connection
between the operator $T$ and the IDE~\eqref{eqn-ide}.

\begin{prp}\label{prp-fixed}
  A function $g\in\cC$
  is a fixed point of the operator $T$
  if and only if 
  $g$ solves the IDE~\eqref{eqn-ide}
  for some $\gamma>0$,
  and then $g$ also satisfies the following properties:
  \begin{substate}
  \item $g\in\breve\cD\`$.
  \item $\area g = \gamma$.
  \item $\stride g = \gamma$.
  \item $g^*$ and $g'$ are continuously differentiable
    on the interval $(0,1]$.
  \item $g''(1)=1$ and $(g^*)'(1) = -\gamma$.
  \end{substate}
\end{prp}
\begin{proof}
  First we assume that $g=Tg\in\cC$
  and set $\alpha:=\area g$.
  Using Proposition~\ref{prp-DIT}(a),
  we conclude that $g\in\breve\cD$
  and then $g\in\breve\cD\`$ by~\ref{rmk-inv}(c)
  and Proposition~\ref{prp-DIT}(d),
  settling assertion~(a).
  Differentiating the equation
  $g(x)=(Tg)(x)=\frac1\alpha\int_x^1g^*$,
  we arrive at $-\alpha g'=g^*$, that is,
  $g$ solves~\eqref{eqn-ide} with $\gamma=\alpha$.
  
  Conversely assume that $g\in\cC$
  (is differentiable and) solves~\eqref{eqn-ide}
  for some $\gamma>0$.
  Integrating~\eqref{eqn-ide}
  while considering~\ref{rmk-inv}(c)
  leads to $g(x) = \frac1\gamma\int_x^1g^*$.
  Plugging $0$ into this,
  yields $1=g(0)=\frac1\gamma\area g$
  by~\ref{rmk-inv}(f),
  thereby showing~(b) and $g=Tg$.
  \begin{subproof}\setcounter{substate}{2}
  \item Plugging $0$ into~\eqref{eqn-ide} 
    and using~\ref{rmk-inv}(c) gives
    $g'(0)=\frac{-g^*(0)}\gamma=\frac{-1}\gamma$,
    hence $\stride g=\gamma$ by~\eqref{eqn-slopes}.
  \item By~(a) and~\ref{rmk-inv}(e),
    we have $g^*\in\cD'$,
    and the assertion follows from~\eqref{eqn-ide}.
  \item Plugging $1$ into the derivative
    of~\eqref{eqn-ide}, yields
    $-\gamma g''(1) = (g^*)'(g(0))
     = \frac1{g'(0)} = -\gamma$
    by the chain rule and~(c),
    thus $g''(1) = 1$ and $(g^*)'(1) = -\gamma$.
    \qedhere
  \end{subproof}
\end{proof}

Next, we want to construct
a complete $T$-invariant subset $\cK$ of~$\breve\cC$.
To this end, we need to bound area and stride
of $Tg$ from below.

\begin{lem}\label{lem-Tgbound}
  Let $g\in\breve\cC$, $0<\alpha\leq\stride g$,
  $\beta := \inf g^{-1}\{0\}$ and $\gamma := \area g$.
  Then
  \begin{substate}
  \item $\alpha\leq2\gamma\leq\beta\leq1$,
    and $\alpha=2\gamma \implies 2\gamma=\beta
    \implies \area Tg=\frac13$,
  \item $(Tg)':[0,1]\to(-\infty,0]$ exists, is continuous,
    strictly increasing and concave,
  \item $Tg\in\breve\cD\`$ with $\stride Tg = \frac\gamma\beta$,
  \item $\frac\beta\gamma(\id-1)\leq(Tg)'
         \leq\frac\alpha\gamma(\id-1)$,
  \item $\area Tg \leq \frac13$,
  \item $\alpha\beta-4\alpha\gamma+4\gamma^2
         \leq 6(\beta-\alpha)\gamma \area Tg$.
  \end{substate}
\end{lem}
\begin{proof}
  \begin{subproof}
  \item From $g\in\breve\cC$ and the definition of $\beta$,
     we infer that $g(x)\leq 1-\frac x\beta$
     for all $x\in[0,\beta]$,
     hence $\frac\alpha2 = \int_0^{\alpha}(1-\frac\id\alpha)
            \leq \area g = \gamma
            \leq \int_0^\beta(1-\frac\id\beta)
            = \frac\beta2$,
     settling the asserted inequality chain.
     From this, we also see that
     $\frac\alpha2=\gamma \iff g\rest{[\alpha,1]}=0
      \implies \beta=\alpha$ and that $\gamma=\frac\beta2
      \implies g^*(y)=\beta(1-y)$ for $y\in(0,1]
      \implies Tg=(1-\id)^2 \implies \area Tg=\frac13$.
  \item By its convexity,
    $g$ is strictly decreasing on $[0,\beta]$.
    Thus $f(x):=g(\beta x)$ for $x\in[0,1]$
    defines a function $f\in\breve\cD$, which satisfies
    $1-\frac\beta\alpha\id \leq f \leq 1-\id$.
    Using~\ref{rmk-inv}(d), (a) and~(f), we infer that
    $f^*\in\breve\cD$,
    \begin{equation}\label{eqn-f*}
      \tfrac\alpha\beta(1-\id) \leq f^* \leq 1-\id
      \quad\text{and}\quad
      \textstyle
      \area f^* = \area f = \frac1\beta\area g
                = \frac\gamma\beta.
    \end{equation}
    Because $\beta f^*(x)=g^*(x)$ for all $x\in(0,1]$,
    we conclude that $Tg=Tf$
    is differentiable with continuous derivative
    \begin{equation}\label{eqn-Tg'}
      (Tg)' = (Tf)' = -\tfrac\beta\gamma f^*,
    \end{equation}
    and the assertions follow.
  \item From~\eqref{eqn-slopes} and~\eqref{eqn-Tg'},
    we infer that
    $\stride Tg = -\frac1{(Tg)'(0)} = \frac\gamma\beta$
    and $(Tg)'(1) = 0$,
    hence $Tg\in\cD\`$,
    while $Tg\in\breve\cD$ holds
    by Proposition~\ref{prp-DIT}(a).
  \item follows from~\eqref{eqn-f*} and~\eqref{eqn-Tg'}.
  \item By (a)--(c)
    and becaus
    $\area(Tg)'=(Tg)(1)-(Tg)(0)=-1=\area(2\id-2)$,
    \[
      s:=\sup\set{0<x<1}{(Tg)'(x)\leq2x-2} \in (0,1]
    \]
    is well-defined,
    $(Tg)'\rest{[0,s]}\leq2\id_{[0,s]}-2$
    and $(Tg)'\rest{[s,1]}\geq2\id_{[s,1]}-2$.
    We conclude that
    $(Tg)(x) \leq 1 + \int_0^x(2\id-2) = (1-x)^2$
    for $x\in[0,s]$ and also
    $(Tg)(x) = -\int_x^1(Tg)' \leq -\int_x^1(2\id-2)
             = (1-x)^2$
    for $x\in[s,1]$.
    Hence, $\area Tg \leq \area(1-\id)^2 = \tfrac13$.
  \item Using (a),
    the asserted inequality is verified directly
    if $\alpha\leq2\gamma=\beta$,
    and we may assume $\alpha<2\gamma<\beta$.
    We conclude that
    $\xi:=\frac{2\gamma-\alpha}{\beta-\alpha}\in(0,1)$
    and define
    $b\colon [0,1]\to\R$ by setting
    \[
      b(x) := \left\{ \begin{array}{ll}
        b_0(x) := 1 -\frac\beta\gamma x
          + \frac{\beta^2-2\alpha(\beta-\gamma)}
                 {2(2\gamma-\alpha)\gamma} x^2
        & \text{for }x\in[0,\xi], \\
        b_1(x) := \frac\alpha{2\gamma} (1-x)^2
        & \text{for }x\in[\xi,1].
      \end{array} \right.
    \]
    It is straightforward to verify that $b\in\breve\cD\`$
    with derivative $b'\colon [0,1]\to\R$ given by
    \[
      b'(x) = \left\{ \begin{array}{ll}
        b_0'(x) = -\frac\beta\gamma
          + \frac{\beta^2-2\alpha(\beta-\gamma)}
                 {(2\gamma-\alpha)\gamma} x
        & \text{for }x\in[0,\xi], \\
        b_1'(x) = -\frac\alpha\gamma(1-x)
        & \text{for }x\in[\xi,1],
      \end{array} \right.
    \]
    which is concave and
    consists of two lines meeting in the point
    $(\xi,\frac\alpha\gamma(\xi-1))$.
    Using (b) and (d), we infer that
    \[
      s:=\inf\set{x\in(0,1]}{(Tg)'(x)\leq b'(x)}
         \in [0,\xi],
    \]
    $(Tg)'\rest{[0,s]}\geq b'\rest{[0,s]}$
    and $(Tg)'\rest{[s,1]}\leq b'\rest{[s,1]}$.
    Thus
    $(Tg)(x) = 1 + \int_0^x(Tg)' \geq 1 + \int_0^x b'
             = b(x)$
    for $x\in[0,s]$ and also
    $(Tg)(x) = -\int_x^1(Tg)' \geq -\int_x^1 b' = b(x)$
    for $x\in[s,1]$,
    hence
    \[
      \area Tg \geq \int_0^1 b
      = \int_0^\xi b_0 + \int_\xi^1 b_1
      = \frac{\alpha\beta-4\alpha\gamma+4\gamma^2}
             {6(\beta-\alpha)\gamma}
    \]
    after a tedious but straightforward calculation.
    \qedhere
  \end{subproof}
\end{proof}

We are now ready to establish the set
\[
  \textstyle
  \cK := \set{g\in\breve\cC}
          {\stride g, \area g\geq\frac15}.
\]
and its properties concerning the operator~$T$,
needed to prove our main theorems.

\begin{thm}\label{thm-K}
The set $\cK$ has the following properties.
	\begin{substate}
	\item $T(\cK)\sseq \cK \cap \breve\cD\`$.
	\item For each $f\in\cE$ there exists $n\in\N$
	  such that $T^nf\in\cK$.
	\item The two metrics $d_\infty$ and $d_1$ are equivalent on $\cK$
	  in the sense that
	  $d_1(f,g) \leq d_\infty(f,g) \leq 5\sqrt{d_1(f,g)}$
	  for all $f,g\in\cK$.
	\item The metric spaces $(\cK,d_\infty)$
	  and $(\cK,d_1)$ are complete.
	\item The restriction $T\rest\cK\colon\cK\to\cK$
	  is continuous when equipping domain and codomain
	  independently with $d_\infty$ or $d_1$.
	\item Every sequence in $(\cK,d_\infty)$ or in $(\cK,d_1)$
	  has a convergent subsequence.
	\end{substate}
\end{thm}
\begin{proof}
	\begin{subproof}
	\item Let $g\in\cK$.
		Then $\gamma:=\area g$, $\stride g\geq\frac15$
		and 	$\beta:=\inf g^{-1}\{0\}\in[\frac25,1]$
		by~\ref{lem-Tgbound}(a).
		With~\ref{lem-Tgbound}(f)
		we infer that $\area Tg \geq u(\gamma)$,
		where the function
		$u\colon (0,\infty)\to\R$ satisfies
	  \[
	    u(x) = \frac{\beta-4x+20x^2}{6(5\beta-1)x}
	    \quad\text{and}\quad
	    u'(x) = \frac{20x^2-\beta}{6(5\beta-1)x^2}
	    \quad\text{for all }x>0.
	  \]
	  We conclude that
	  $\area Tg \geq u\bigl(\frac\beta{2\sqrt5}\bigr)
     = \frac23\cdot\frac{\sqrt5-1}{5\beta-1}
     \geq \frac{\sqrt5-1}6
     > \frac15$.
    Moreover,
    $Tg\in\breve\cD\`\sseq\breve\cC$ and
    $\stride Tg = \frac\gamma\beta \geq \gamma \geq \frac15$
		by~\ref{lem-Tgbound}(c).
	  In total we have shown $Tg\in\cK \cap \breve\cD\`$.
	\item Let $f\in\cE$, and set $f_n:=T^{2+n}f$,
	  $\gamma_n:=\area f_n$ and $\alpha_n := \stride f_n$
    for $n\in\N_0$.
    For $0<\alpha<1$,
    define $\thet(\alpha):=\frac32(\alpha+\sqrt\alpha)$
    and the function
	  $u_\alpha\colon (0,\infty)\to(0,\infty)$
	  given by
	  \[
      u_\alpha(\gamma)
      := \frac{\alpha-4\alpha\gamma+4\gamma^2}
              {6(1-\alpha)\gamma},
      \quad\text{thus}\quad
	    u_\alpha'(\gamma)
      = \frac{4\gamma^2-\alpha}{6(1-\alpha)\gamma^2}
	    \quad\text{for all }\gamma>0,
    \]
	  which therefore satisfies
    $u_\alpha(\gamma)
    \geq u_\alpha\bigl(\frac{\sqrt\alpha}2\bigr)
    = \frac\alpha{\thet(\alpha)}$
	  for all $\gamma>0$.
    With Proposition~\ref{prp-DIT}(a) and~(b)
    and Lemma~\ref{lem-Tgbound}(c) and~(f),
    we obtain
    \[
      f_n\in\breve\cD,
      \quad
      \alpha_{n+1} = \gamma_n
      \quad\text{and}\quad
      \gamma_{n+1} \geq u_{\alpha_n}(\gamma_n)
                   \geq \tfrac{\alpha_n}{\thet(\alpha_n)}
      \quad\text{for all }n\in\N_0.
    \]
    Because
    \[
      0<\alpha\leq\tfrac15 \implies
      \thet(\alpha) \leq \thet\bigl(\tfrac15\bigr) < 1
      \quad\text{and}\quad
      \tfrac15\leq\alpha<1 \implies
      \tfrac\alpha{\thet(\alpha)} \geq
      \tfrac{\sqrt5-1}6 > \tfrac15,
    \]
    there consequently exists $n\in\N_0$
    with $\alpha_n,\gamma_n\geq\frac15$,
    hence $f_n\in\cK$.
	\item Let $f,g\in\cK$. The estimate
	  $\area|f-g|\leq\sup_{x\in[0,1]}|f(x)-g(x)|$
	  settles the left inequality.
	  As for the right one, we may assume that
	  $\delta:=d_\infty(f,g) = f(x_0)-g(x_0)$
    for some $x_0\in[0,1]$.
    From $\stride f,\stride g\geq\frac15$ and $f,g\in\breve\cC$
    we infer that
    \[
      f(x)-g(x)\geq \delta-5|x-x_0|
      \quad\text{for all } x\in[0,1],
    \]
    hence
    $0\leq a:=x_0-\frac\delta5\leq b:=x_0+\frac\delta5\leq1$
    and
    $d_1(f,g) \geq \int_a^b(f-g) \geq \frac\delta2(b-a)
     = \frac{\delta^2}5$.
	\item Recall that $\cC^0[0,1]$,
	  the $\R$-vector space of continuous functions
	  on the interval $[0,1]$,
	  is complete with respect to the $\sup$-norm.
    Therefore, each Cauchy sequence $(g_n)_{n\in\N}$
    in~$(\cK,d_\infty)$ converges to some function $g\in\cC^0[0,1]$
	  satisfying $g(0)=1$ and $g(1)=0$.
	  Clearly, $g$ is again decreasing and convex,
	  and both inequalities
	  $g\geq1-5\id$ and $\area g\geq\frac15$ hold.
	  Hence~$g\in\cK$, showing that $(\cK,d_\infty)$ is complete.
	  The completeness of $(\cK,d_1)$ follows with~(c).
	\item Let $(g_n)_{n\in\N}$ be a sequence in $\cK$
	  converging to $g\in\cK$
	  with respect to $d_\infty$ or $d_1$.
	  Setting $\hat g:=\area g\cdot Tg$
	  and $\hat g_n:=\area g_n\cdot Tg_n$,
	  Remark~\ref{rmk-inv}(g) yields
	  \[\textstyle
	    |\hat g(x)-\hat g_n(x)|
	    = \bigl|\int_x^1(g^*-g_n^*)\bigr|
	    \leq \area|g^*-g_n^*|
	       = d_1(g,g_n) \leq d_\infty(g,g_n)
	  \]
	  for all $x\in[0,1]$, implying
	  $\lim_{n\to\infty}d_\infty(\hat g,\hat g_n)=0$ and
	  $\lim_{n\to\infty}\area g_n = \area g$
	  because $(Tg)(0)=1=(Tg_n)(0)$.
	  We conclude that
	  $d_1(Tg,Tg_n)\leq d_\infty(Tg,Tg_n)\to0$ as $n\to\infty$.
	\item Let $g_n\in\cK$ for $n\in\N$.
	  Because $g_n$ is convex and $\stride g_n\geq\frac15$,
	  we conclude that $|g_n(x_1)-g_n(x_2)|\leq5|x_2-x_1|$
	  for all $x_1,x_2\in[0,1]$, $n\in\N$.
	  Therefore the sequence $(g_n)_{n\in\N}$ is
	  uniformly equicontinuous,
	  and as it is also uniformly bounded,
	  the Arzelà-Ascoli theorem guarantees that
	  it has a convergent subsequence in $(\cK,d_\infty)$.
	  By~(d) its limit lies in $\cK$,
	  and by (c) the same subsequence
	  converges in $(\cK,d_1)$ to the same limit.
	\qedhere
	\end{subproof}
\end{proof}

Henceforth, when speaking about convergence in $\cK$,
we mean convergence in $(\cK,d_\infty)$,
i.e.\ uniform convergence, {\em and}, by~\ref{thm-K}(c) equivalently,
convergence in $(\cK,d_1)$.

\section{Crossing number}

The crucial proofs of~\ref{lem-Tgbound}(e) and (f)
rest on the fact
that $(Tg)'$ intersects another derivative at most once.
More generally, if $f,g\in\cD$, $g\leq f\neq g$ and
$(Tg)'-(Tf)'= \frac{f^*}{{\int}f} - \frac{g^*}{{\int}g}$
changes its sign only once (from $-$ to $+$ in this case),
then we will have $Tg\leq Tf$.
To propagate this reasoning to the next iteration step,
we would require the difference of $(Tf)^*$ and $(Tg)^*$,
after somehow stretching them vertically, to also
change sign at most once. But a vertical stretching of,
say, $(Tg)^*$ corresponds to a horizontal stretching of $Tg$
and thus of $(Tg)'= -\frac{g^*}{{\int}g}$,
which again corresponds to a
vertical {\em and} horizontal stretching of~$g$.
Because it is hard to tell the stretching factors in advance,
we will consider the difference of $f$ and $g$ after
arbitrary horizontal {\em and} vertical stretching.

As a first step, we want to count
how often a given continuous function
$\Delta\colon[a,b]\to \R$
defined on a bounded, closed interval $[a,b]$
changes sign.
To this end, we call a closed subinterval
$[c,d]\sseq[a,b]$ with $a<c\leq d<b$
and image $\Delta([c,d])=\{0\}$
a \df{sign switch} of $\Delta$ if there exists
$\delta\in(0,\min\{c-a,b-d\}]$ such that
$\Delta(c-x)\cdot\Delta(d+x)<0$
for all $x\in(0,\delta]$.
By $\cX\Delta$ we denote the set
of all sign switches of $\Delta$
and by $\chi\Delta := \#\cX\Delta  \in\N_0\cup\{\infty\}$
their number.

\begin{rmk}\label{rmk-Delta}
  Let $k\in\N_0$, $a,b,a',b'\in\R$ with $a<b$ and $a'<b'$.
  Let $u\colon[a',b']\to[a,b]$ and $\Delta\colon[a,b]\to\R$
  be continuous functions, $u$ bijective.
  The following statements hold.
  \begin{substate}
  \item $\cX(c\Delta)=\cX\Delta$ for every $c\in\R\sm\{0\}$.
  \item $\chi(\Delta\circ u)=\chi\Delta$.
  \item $\chi\Delta\geq k$ if and only if there exist
    $a\leq x_0<\dotsm<x_k \leq b$
    such that $\Delta(x_{i-1})\cdot\Delta(x_i)<0$
    for $i\in\enum1k$.
  \item If $\Delta(a)\cdot\Delta(b)>0$,
    then $\chi\Delta$  is even or~$\infty$.
  \item If $\Delta(a)\cdot\Delta(b)<0$,
    then $\chi\Delta$ is odd or~$\infty$.
  \item Suppose that $\Delta$
    is continuously differentiable.
    Then $\chi\Delta'\geq\chi\Delta-1$.
    If  $\Delta(a)\cdot\Delta'(a)>0$ in addition,
    then $\chi\Delta'\geq\chi\Delta$.
  \item Let $\Delta_n\colon[a,b]\to\R$ for $n\in\N$
    be continuous functions such that
    $(\Delta_n)_{n\in\N}$ converges to~$\Delta$.
    Then $\chi\Delta\leq\sup_{n\in\N}\chi\Delta_n$.
  \end{substate}
\end{rmk}
\begin{proof}
   {\bf(a)}--{\bf(c)} are immediate
   from the definition of sign switches.
  \begin{subproof}
  \setcounter{substate}3
  \item and {\bf(e)} follow from (c).
  \setcounter{substate}5
  \item Suppose that $k:=\chi\Delta\in\N_0$.
	  Then there are $a\leq x_0<\dotsm<x_k \leq b$ as in~(c).
	  By the mean value theorem, we
	  can find $y_i\in[x_{i-1},x_i]$ with
	  $\Delta'(y_i)\cdot\Delta(x_i)>0$ for $i\in\enum1k$.
	  This shows $\chi\Delta'\geq k-1$
	  according to~(c).
	  
    If $\Delta(a)\cdot\Delta'(a) > 0$,
    then we can find $y_0\in[a,x_0]$
    with $\Delta'(y_0)\cdot\Delta(x_0)>0$;
    hence $\chi\Delta'\geq k$ by~(c) again.
  \item Let $\chi\Delta\geq k$
    with $\thru x0k$ as in (c).
    By assumption, we can choose $n\in\N$ such that
    $|\Delta_n(x_i)-\Delta(x_i)| < |\Delta(x_i)|$,
    hence $\Delta_n(x_i)\cdot\Delta(x_i)>0$
    for $i\in\enum0k$.
    This implies
    $\Delta_n(x_{i-1})\cdot\Delta_n(x_i)<0$
    for $i\in\enum1k$ which, by~(c), is equivalent to
    $\chi\Delta_n\geq k$.
    \qedhere
  \end{subproof}
\end{proof}

Given two functions $f,g\in\cD$ and $a,b>0$,
we introduce the function
\[
  f\mul a := f\circ(a\id) \colon [0,\tfrac1a]\to[0,1],
  \quad x\mapsto f(ax)
\]
obtained by stretching $f$ horizontally with the factor $\frac1a$,
and consider the continuous function
$f\mul a-bg \colon [0,\min\{1,\tfrac1a\}]\to\R$.
The next lemma tells us how its number of sign switches
behaves under swapping $f$ with $g$
and under the operators ${}^*$ and~$I$.

\begin{lem}\label{lem-chi}
  Let $a,b>0$ and  $f,g\in\cD$.
  The following statements hold.
  \begin{substate}
  \item $\chi(f\mul a-bg)
       = \chi(g\mul\frac1a-\frac1bf)$.
  \item $\chi(f\mul a-bg)
       = \chi(g^*\mul\tfrac1b-\tfrac1a\cdot f^*)$.
  \item  Let $\hat\Delta := If\mul a - bIg$,
    $b':= \frac{b{\int} f}{a{\int} g}$
    and $\Delta := f\mul a - b'g$.
    Then $\chi\hat\Delta \leq 1+\chi\Delta$.\\
    If $b<1<b'$ or $b'<1<b$,
    then $\chi\hat\Delta \leq \chi\Delta$.
  \item If either $a,b<1$ or $a,b>1$,
    then $\chi(f\mul a-bg)$ is even or $\infty$.
  \item If $a<1<b$ or $b<1<a$, then $\chi(f\mul a-bg)$
    is odd or $\infty$.
  \item If $a,b\leq1$ and $g\leq f$,
    then $\chi(f\mul a-bg)=0$.
  \end{substate}
\end{lem}

\begin{proof}
  \begin{subproof}
  \item Let $a':=\min\{1,\frac1a\}$,
    $\Delta:=f\mul a-bg$
    and $\ti\Delta:=g\mul\frac1a-\frac1bf$.
    Then
    \[
      b\ti\Delta(ax)
      = b\cdot\bigl(g(\tfrac{ax}{a})-\tfrac1bf(ax)\bigr)
      = bg(x)-f(ax)
      = -\Delta(x)
    \]
    for all $x\in[0,a']$.
    Hence, $\chi\Delta=\chi\ti\Delta$
    by~\ref{rmk-Delta}(a) and (b).
  \item Let  $a':=\min\{1,\frac1a\}$, $b':=\min\{1,b\}$,
    $\Delta:=f\mul a-bg\colon[0,a']\to\R$
    and
    \[
      \ti\Delta :=
      g^*\mul\tfrac1b-\tfrac1af^*
      \colon [0,b']\to\R.
    \]
    Because the function
    $u\colon[0,a']\to[0,b']$,
    $x\mapsto \min\{f(ax),bg(x)\}$
    is bijective by~\ref{rmk-inv}(c)
    and $\cX(\ti\Delta\circ u)=\cX\Delta$,
    the assertion follows from~\ref{rmk-Delta}(b).
  \item According to its definition,
	  $\hat\Delta$ is differentiable
	  with continuous derivative
	  \[
	    \hat\Delta'
	    = a(If)'\mul a - b(Ig)'
	    = -\tfrac{a}{{\int}f}\Delta,
	  \]
	  hence
	  $\chi\Delta = \chi\hat\Delta' \geq \chi\hat\Delta - 1$
	  by~\ref{rmk-Delta}(a) and~(f).
    If $b<1<b'$ or $b'<1<b$, then
    $\hat\Delta(0)\cdot\hat\Delta'(0)
     = (1-b)\cdot\frac{a}{{\int}f}\cdot(b'-1) > 0$,
    hence
	  $\chi\Delta = \chi\hat\Delta' \geq \chi\hat\Delta$,
	  again by~\ref{rmk-Delta}(a) and~(f).
  \item Set $\Delta:=f\mul a-bg$.
    If $a,b<1$, then $\Delta(0)=f(0)-bg(0)=1-b>0$
    and $\Delta(1)=f(a)>0$.
    If $a,b>1$, then $\Delta(0)=f(0)-bg(0)=1-b<0$
    and $\Delta(\frac1a)=-g(\frac1a)<0$.
    In both cases,
    the assertion follows from~\ref{rmk-Delta}(d).
  \item Set $\Delta:=f\mul a-bg$.
    If $a<1<b$, then $\Delta(0)=f(0)-bg(0)=1-b<0$
    and $\Delta(1)=f(a)>0$.
    If $b<1<a$, then $\Delta(0)=f(0)-bg(0)=1-b>0$
    and $\Delta(\frac1a)=-g(\frac1a)<0$.
    In both cases,
    the assertion follows from~\ref{rmk-Delta}(e).
  \item From $f,g\in\cD$, $a,b\leq1$ and $g\leq f$, we conclude
     $f(ax)\geq f(x)$ and $bg(x)\leq g(x)$, hence
     $(f\mul a-bg)(x) = f(ax)-bg(x) \geq f(x)-g(x) \geq 0$
     for all $x\in[0,1]$, so that $\chi(f\mul a-bg)=0$.
  \qedhere
  \end{subproof}
\end{proof}

\comment{
\begin{lem}
  For all $f,g\in\cD$,
  the map $\chi_{f,g}\colon(0,\infty)^2\to[0,\infty)$,
  $(a,b)\mapsto\frac{1}{1+\chi(f\mul a-bg)}$ is continuous.
\end{lem}
}

Given $f,g\in\cD$,
we define the \df{crossing number}
\[
  \chi(f,g)
   := \sup\set{\chi(f\mul a-bg)}{a,b>0}
  \in\N_0\cup\{\infty\}
\]
\df{of $f$ with $g$}.
We write $g\lhd f$ or $f\rhd g$
and say that $f$ \df{dominates}~$g$
if $\chi(f,g)=2$ and $g\leq f$.

\begin{lem}\label{lem-dom}
  For $f,g\in\cD$,
  the following statements hold.
  \begin{substate}
  \item $\chi(f,g)=\chi(g,f)\geq1$.
  \item $g\lhd f \iff g^*\lhd f^*$.
  \item $\chi(f,f)$ is odd or $\infty$.
  \item $\chi(f,g)\leq2 \implies \chi(f-g)=0
         \iff f\leq g$ or $g\leq f \implies f=g$ or $\chi(f,g)\geq2$. 
  \item $\chi(f,g)=2 \iff f\neq g$
    and either $f\lhd g$ or $g\lhd f$.
  \item If $g \lhd f$ and $0<\min\{a,b\}\leq1$,
    then $\chi(f\mul a-bg) \leq 1$.
  \end{substate}
\end{lem}
\begin{proof}
  \begin{subproof}
  \item follows from~\ref{lem-chi}(a) and~(e).
  \item follows from~\ref{lem-chi}(b) and~\ref{rmk-inv}(a).
  \item If $a,b\leq1$ or $a,b\geq1$,
    then $\chi(f\mul a - bf)=0$ by~\ref{lem-chi}(a) and~(f).
    Therefore, the assertion follows with~\ref{lem-chi}(e).
  \item We prove the first implication by contraposition.
    Suppose that $\chi(f-g)>0$,
    that is, neither $f\leq g$ nor $g\leq f$ holds.
    Then we may assume w.l.o.g.
    that there are $0<x_1<x_2<1$
    with $f(x_1)>g(x_1)$ and $f(x_2)<g(x_2)$,
    and $\Delta:=f\mul(1-\eps)-(1+\eps)g$
    by continuity still satisfies
    $\Delta(x_1)>0$ and $\Delta(x_2)<0$ 
    for sufficiently small $\eps>0$.
    But then $\Delta(0)=-\eps<0$
    and $\Delta(1)=f(1-\eps)>0$,
    so that $\chi\Delta\geq3$~by~\ref{rmk-Delta}(c).
    The last implication is obvious.
  \item follows from~(c), (d) and the definition of $\lhd$.
  \item Set $\Delta = f\mul a-bg$, then $\chi\Delta\leq2$.
    The case $a,b\in(0,1]$ is covered by~\ref{lem-chi}(f).
    If $0<a\leq1<b$, then $\Delta(0)=b-1>0$,
    and $\chi(\Delta)=2$ would imply $\Delta(x)<0$ for some $x\in(0,1)$,
    where we may assume $a<1$ by continuity,
    in violation of~\ref{lem-chi}(e).
    If $0<b\leq1<a$, then $\Delta(\frac1a)=-bg(\frac1a)<0$,
    and $\chi(\Delta)=2$ would imply $\Delta(x)>0$
    for some $x\in(0,\frac1a)$,
    where we may assume $b<1$ by continuity,
    violating~\ref{lem-chi}(e) again.
    \qedhere
  \end{subproof}
\end{proof}

Now we can show that domination in $\cD$
is preserved by the operators $I$ and $T$.

\begin{thm}\label{thm-dom}
  Let $f,g\in\cD$ such that $g \lhd f$.
  Then $Ig \lhd If$ and $Tg \lhd Tf$.
\end{thm}
\begin{proof}
  By continuity and~\ref{prp-DIT}(d), we have
  \begin{equation}\label{eqn-area}
    \stride Ig=\area g < \area f = \stride If.
  \end{equation}
  Let $a,b>0$ and
  $b':= \frac{b{\int}f}{a{\int}g}$.
  Then $k' := \chi(f\mul a - b'g)\leq2$
  and $k := \chi(If\mul a - bIg)\leq1+k'$
  according to~\ref{lem-chi}(c).
  Because of~\ref{lem-chi}~(d),
  we need only consider the two cases $a\leq1\leq b$
  and $b<1\leq a$, in order to show that $k\leq2$.
  \begin{mycases}
  \item If $a\leq1\leq b$,
	  then $b'>1$ by~\eqref{eqn-area}, hence
	  $k \leq 1+k' \leq 2$ by~\ref{lem-dom}(f).
  \item Now suppose that $b<1\leq a$.
	  If $b'\leq1$, then $k\leq1+k'\leq2$
	  by~\ref{lem-dom}(f),
    and if $b'>1$, then
    $k \leq k'\leq2$ by~\ref{lem-chi}(c).
  \end{mycases}
  In view of~\ref{lem-dom}(a) and~\ref{lem-chi}(d)
  we have altogether proved $\chi(If,Ig)=2$,
  hence $Ig\lhd If$ by~\eqref{eqn-area}.
  Because $f^*,g^*\in\cD$ by~\ref{rmk-inv}
  and $g^*\lhd f^*$ by~\ref{lem-dom}(b),
  we also obtain $Tg=Ig^*\lhd If^*=Tf$.
\end{proof}

We conclude this section by observing
that domination is also preserved under taking limits.

\begin{prp}\label{prp-fglim}
  Let $f,g,f_n,g_n\in\cD$ for $n\in\N$
  with $f=\lim_{n\to\infty}f_n$ and
  $g=\lim_{n\to\infty}g_n$.
  The following statements hold.
  \begin{substate}
	\item $\chi(f,g)\leq\sup_{n\in\N}\chi(f_n,g_n)$,
	\item If $g_n\lhd f_n$ for all $n\in\N$,
	  then either $f=g$ or $g\lhd f$.  
  \end{substate}
\end{prp}
\begin{proof}
  \begin{subproof}
  \item Let $a,b>0$, and set 
    $\Delta_n := f_n\mul a-bg_n$ for $n\in\N$.
    Then $(\Delta_n)_{n\in\N}$ converges to 
    $\Delta := f\mul a-bg$,
    hence $\chi\Delta\leq\sup_{n\in\N}\chi\Delta_n$
    due to~\ref{rmk-Delta}(g).
  \item From the assumption and~(a), we infer that
    $g\leq f$ and $\chi(f,g)\leq2$.
    With~\ref{lem-dom}(a), (d) and~(e),
    we conclude that either $f=g$ or $g\lhd f$.
    \qedhere
  \end{subproof}
\end{proof}

\section{Existence, uniqueness and global convergence}

We will now construct a unit stribola,
that is, a solution to the IDE~\eqref{eqn-ide}.
To this end, we define the \df{canonical stribolic iterates}
and their \df{areas},
\begin{equation}\label{eqn-hn}
  h_n:= T^n1_{[0,1]} \in \cE,
  \quad\kappa_n := \area h_n
  \quad\text{for }n\in\N_0.
\end{equation}
In particular, we have
$h_0(x) = 1 = h_0^*(x)$,
$h_1(x)=(Ih_0^*)(x) = 1-x = h_1^*(x)$,
$h_2(x)=(Ih_1^*)(x)=(1-x)^2$, $h_2^*(x)=1-\sqrt x$,
$h_3(x)=(Ih_2^*)(x) = 1-3x+2x^{\frac32}$
for $x\in[0,1]$, and 
$(\kappa_0,\kappa_1,\kappa_2,\kappa_3,\dotsc)
= (1,\frac12,\frac13,\frac3{10},\dotsc)$.

\begin{figure}[h]
\centering
\includegraphics[width=0.6\textwidth]{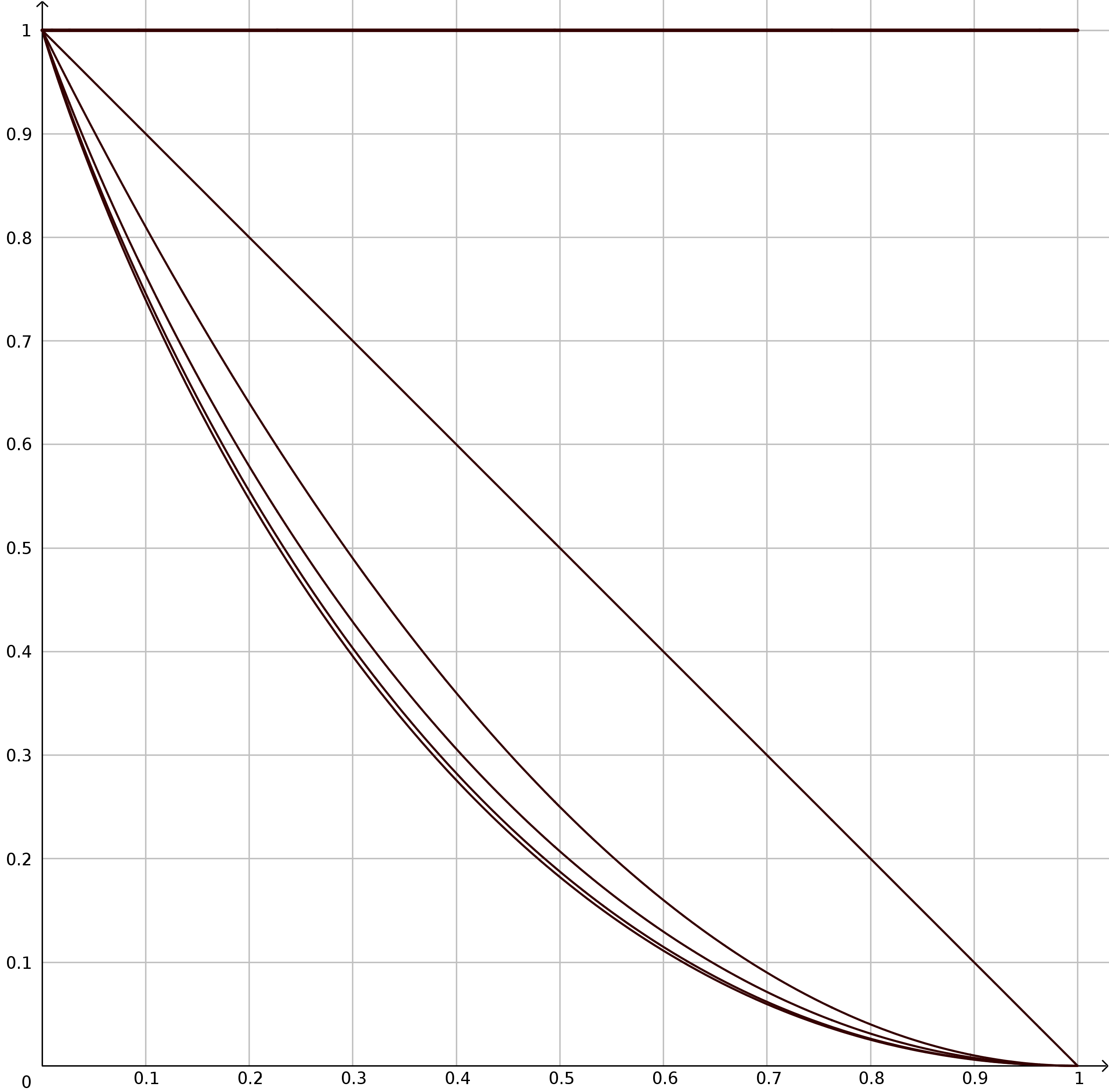}
\caption{Graphs of $\thru h05$}
\label{fig-h012345}
\end{figure}

Repeated application of Theorem~\ref{thm-dom} will show
that these iterates descend to a unit stribola.

\begin{thm}\label{thm-descend}
The sequence $(h_n)_{n\in\N_0}$ from~\eqref{eqn-hn}
satisfies $h_{n+1}\lhd h_n\in\cK$ for all $n\in\N$
and converges to a function $h\in\cK$.
Moreover, $h = Th$ solves~\eqref{eqn-ide}
with $\gamma=\area h$.
\end{thm}
\begin{proof}
  Because $h_1\in\cK$ is linear,
  it dominates the non-linear, convex function $h_2\in\cK$.
  With Theorems~\ref{thm-dom}
  and~\ref{thm-K}(a),
  we inductively conclude that
  \begin{equation}\label{eqn-descend}
    h_{n+1} \lhd h_n \in \cK
    \quad\text{for all } n\in\N.
  \end{equation}
  By~\ref{thm-K}(f),
  there are positive integers $n_1<n_2<\dotsm$
  and $h\in\cK$ such that
  $\lim_{k\to\infty}h_{n_k} = h$,
  which implies $\lim_{n\to\infty}h_n = h$
  because of~\eqref{eqn-descend}.
  Finally, using~\ref{thm-K}(e), \eqref{eqn-hn}
  and Proposition~\ref{prp-fixed}, we see that
  \[
    Th = T\Bigl(\lim_{n\to\infty}h_n\Bigr)
            = \lim_{n\to\infty}Th_n
            = \lim_{n\to\infty}h_{n+1} = h
  \]
  solves~\eqref{eqn-ide} with $\gamma=\area h$.
\end{proof}

Until the end of this paper,
we shall focus on the unit stribola
$h := \lim_{n\to\infty}h_n \in\breve\cD\`$
and its area
\begin{equation}\label{eqn-kappa}
  \kappa := \area h = -\frac{1}{h'(0)} = \stride h
  = \lim_{n\to\infty}\kappa_n = \inf_{n\in\N_0}\kappa_n,
\end{equation}
as established in Theorem~\ref{thm-descend}.

\begin{cor}\label{cor-hdom}
  For $m,n\in\N$ with $m<n$, we have
  \begin{substate}
  \item $h_n\lhd h_m$,
  \item $h\lhd h_m$,
  \item $\chi(h_m,h_m)=1$,
  \item $\chi(h,h)=1$.
  \end{substate}
\end{cor}
\begin{proof}
  \begin{subproof}
  \item Clearly, $h_n\lhd h_1$ because $h_1$ is linear,
    $h_n\neq h_1$ and $h_n$ is convex.
    The assertion follows
    with Theorem~\ref{thm-dom} by induction on $m$.
  \item holds according to~(a) and Proposition~\ref{prp-fglim}(b).
  \item We proceed by induction on $m$.
    Obviously, $\chi(h_1,h_1)=1$.
    Assume that $\chi(h_m,h_m)=1$.
    Because $h_{m+1} = Ih_m^*$, we conclude
    that $\chi(h_{m+1},h_{m+1})\leq2$
    using~\ref{lem-chi}(b) and~(c),
    hence $\chi(h_{m+1},h_{m+1})=1$ by~\ref{lem-dom}(c).
  \item follows from~(c) with~\ref{prp-fglim}(a)
    and~\ref{lem-dom}(a).\qedhere
  \end{subproof}
\end{proof}

As an aside, we observe that the function
\begin{equation}\label{eqn-tih}
  \ti h:=\tfrac1\kappa h\mul\kappa\,\colon
  [0,\tfrac1\kappa]\to[0,\tfrac1\kappa]
  \quad\text{satisfies}\quad -\ti h'=\ti h^*,
\end{equation}
that is,
$\ti h$ becomes its own derivative when rotated clockwise
about the origin by $90^\circ$.
A function defined on $[0,a]$ for some $a>0$
and satisfying the above IDE
shall be called a \df{standard stribola}.

Our next goals are to prove
that $h$ is the only unit stribola
and that $(T^nf)_{n\in\N}$ converges
to $h$ for every $f\in\cE$ (global convergence).
We begin with a result frequently used in the sequel
that involves the stride,
see~\eqref{eqn-stride} and~\eqref{eqn-slopes}.

\begin{stridelem}\label{lem-stride}
  Let $f,g\in\breve\cD$ satisfy $g\lhd f$ and $\stride g>0$.
  Then
  \begin{equation}\label{eqn-mulstride}
    f\mul \stride f \leq g\mul \stride g
  \end{equation}   
  and
  $\stride g\cdot\area f < \stride f\cdot\area g$.
\end{stridelem}
\begin{proof}
  We have $0<\stride g\leq\stride f$ by assumption.
  Suppose~\eqref{eqn-mulstride} were wrong.
  Then $f(\frac{\stride f}{\stride g}x_2)>g(x_2)$
  for some $x_2\in(0,\frac{\stride g}{\stride f})$.
  By continuity, there is $a>\frac{\stride f}{\stride g}$
  such that $f(ax_2)>g(x_2)$ still holds.
  But because $\stride(f\mul a)=\frac{\stride f}a<\stride g$,
  we can find $x_1\in(0,x_2)$ with $f(ax_1)<g(x_1)$.
  Again by continuity, there is $b\in(0,1)$ such that
  $\Delta := f\mul a - bg$ still satisfies
  $\Delta(x_1)<0$ and $\Delta(x_2)>0$.
  Since also $\Delta(0)=1-b>0$
  and $\Delta(\frac1a)=-bg(\frac1a)<0$,
  we would have $\chi(f,g)\geq3$ by~\ref{rmk-Delta}(c),
  contradicting the assumption $g\lhd f$.
  Therefore~\eqref{eqn-mulstride} holds.
  Because $f\mul \stride f = g\mul \stride g$
  would imply $f=g$, we conclude that
  $\stride g\cdot\area f < \stride f\cdot\area g$.
\end{proof}

\begin{cor}\label{cor-stride}
  If the unit stribola $h$ dominates $g\in\breve\cD$,
  then $\stride g<\area g = \stride Tg < \kappa$.
\end{cor}
\begin{proof}
  We may assume that $\stride g>0$
  and obtain
  $\kappa\cdot\stride g = \stride g\cdot\area h <
   \stride h\cdot\area g = \kappa\area g$
  from~\ref{lem-stride},
  hence
  $\stride g<\area g = \stride Tg$
  by~\ref{prp-DIT}(d) and~\ref{rmk-inv}(f).
  Moreover, $g\lhd h$ implies $\area g<\area h=\kappa$.
\end{proof}

\begin{cor}\label{cor-kappan}
  The sequence
  $\bigl(\frac{\kappa_n}{\kappa_{n-1}}\bigr)_{n\in\N}$
  is strictly increasing and converges to $1$.
\end{cor}
\begin{proof}
  Let $n\in\N$.
  Applying the Stride Lemma~\ref{lem-stride}
  to $h_{n+1}\lhd h_n$ yields
  \[
  \kappa_n^2 =
   \stride h_{n+1}\cdot\area h_n < \stride h_n\cdot\area h_{n+1}
   = \kappa_{n-1}\kappa_{n+1},
  \]
  hence
  $\frac{\kappa_n}{\kappa_{n-1}}<\frac{\kappa_{n+1}}{\kappa_n}$.
  The convergence to $1$ holds because of~\eqref{eqn-kappa}.
\end{proof}

We shall use the following lemma as a tool to estimate areas from below.

\begin{hammocklem}\label{lem-hammock}
  Let $\ti f,f,g,\in\breve\cD$ and $c:=\area f / \area g$
  satisfy $g\lhd f$, $\ti f\leq f$
  and $\ti f\leq \ti g := g\mul\frac1c$.
  Then
  $\area Tf - \area Tg < 1 - \area\ti f / \area f$.
\end{hammocklem}
\begin{proof}
  By assumption, we have $c>1$ and
  \begin{gather}\label{eqn-tig}
      c\area g^*
      = c\area g
      = \textstyle\int_0^c\ti g
      = \area f = \area f^*.  
  \end{gather}
  Because $\chi(f,g)\leq2$,
  there is $x_1\in(0,1)$ such that
  $f\rest{[0,x_1]}\geq \ti g\rest{[0,x_1]}$
  and $f\rest{[x_1,1]}\leq \ti g\rest{[x_1,1]}$.
  Therefore $y_1:=\ti g(x_1)=f(x_1)$,
  $\delta:=\int_{y_1}^1(f^*-cg^*) = \int_0^{x_1}(f-\ti g)$
  and $\Delta := Tf-Tg$
  satisfy $\Delta(0)=0=\Delta(1)$ and
  \[
    0 \leq \area f \cdot \Delta
      \leq \area f \cdot \Delta(y_1) = \delta
      < \textstyle\int_0^1(f-\ti f)
      = \area f - \area\ti f
  \]
  by~\eqref{eqn-tig} and the assumptions.
  The asserted inequality follows
  upon integration.
\end{proof}

\begin{cor}\label{cor-hammock}
  Let $n\in\N$ and $g\in\breve\cD$ satisfy $g\lhd h,h_n$ and
  $\kappa\area g \leq \kappa_n\cdot\stride g$.
  Then \mbox{$\kappa_{n+1}-1+\frac{\kappa}{\kappa_n} < \area Tg$}.
\end{cor}
\begin{proof}
  From $c := \frac{{\int}h_n}{{\int}g}
           = \frac{\kappa_n}{{\int}g}
        \geq \frac{\kappa}{\stride g}           
           = \frac{\stride h}{\stride g}$
  we conclude that $h \leq \ti g := g\mul\frac1c$
  by the Stride Lemma~\ref{lem-stride}
  and can apply Hammock Lemma~\ref{lem-hammock}
  with $\ti f=h$ and $f=h_n$.
\end{proof}

Due to Theorem~\ref{thm-K},
for any $f\in\cE$, the limit set
\[
  \cL(f) := \set{\lim_{k\to\infty}T^{m_k}f}
          {0<m_1<m_2<\dotsm \text{ and }
          (T^{m_k}f)_{k\in\N} \text{ is Cauchy}}
\]
of the sequence $(T^nf)_{n\in\N_0}$
is well-defined and satisfies
$\emptyset\neq\cL(f)\sseq\cK \cap \breve\cD\`$.

\begin{lem}\label{lem-sublim}
  Let $f\in\cE$ and $g\in\cL(f)$. Then
  \begin{substate}
  \item $T^mg\in\cL(f)$ for all $m\in\N$,
  \item $g\lhd h_n$ for all $n\in\N$,
  \item $g=h$ or $g\lhd h$,
  \item $h\in\cL(f)$ implies $\lim_{n\to\infty}T^nf=h$.
  \end{substate}
\end{lem}
\begin{proof}
  By Theorem~\ref{thm-K}(a)--(b),
  we may assume that $h_1\neq f\in\cK \cap\breve\cD\`$
  and conclude that
  \begin{equation}\label{eqn-fnhn}
    f_n:=T^{n-1}f\lhd\thru h1n  
  \end{equation}
  for all $n\in\N$
  using Theorem~\ref{thm-dom} with induction.
  By assumption, there are positive integers
  $m_1<m_2<\dotsm$ such that
  $g=\lim_{k\to\infty}f_{m_k}$,
  implying~(b) and~(c)
  by~\eqref{eqn-fnhn} and Proposition~\ref{prp-fglim}(b).
  Due to~\ref{thm-K}(e),
  $T^m$ is continuous on $\cK$, hence
  $T^mg = \lim_{k\to\infty}T^mf_{m_k}
        =\lim_{k\to\infty}f_{m_k+m} \in \cL(f)$
  for all $m\in\N$, proving~(a).
  As for~(d), let $\eps>0$.
  Since $h\in\cL(f)$ and $(\kappa_n)_{n\in\N_0}$
  is decreasing with limit $\kappa$,
  we can choose $m\in\N$ such that
  \[
    \area f_m > \kappa-\eps
    \quad\text{and}\quad
    \kappa_m < \frac{\kappa-\eps}{\kappa-2\eps}\kappa.
  \]
  Using~\eqref{eqn-fnhn} and
  applying the Stride Lemma~\ref{lem-stride}
  to $f_{m+1}\lhd h_{m+1},\dotsc,f_n\lhd h_n$
  yields
  \[
    \area f_n > \frac{\kappa_n}{\kappa_m}\area f_m
             > \frac{\kappa}{\kappa_m}\area f_m
             > \frac{\kappa-2\eps}{\kappa-\eps}\area f_m
             > \kappa-2\epsilon
  \]
  for all $n>m$.
  We conclude that $\lim_{n\to\infty}f_n = h$
  using~\eqref{eqn-fnhn} again.
\end{proof}

Although the strong uniqueness established in the next theorem
would follow from the global convergence,
we feel like proving it directly
because it drops out easily
from a small subset of our previous results.

\begin{thm}\label{thm-unique}
  Suppose that $r\in\N$ and $g\in\cE$ satisfy $T^rg=g$.
  Then $g=h$.
\end{thm}
\begin{proof}
  The assumptions and \ref{thm-K}(b) imply
  $\cL(g)=\{g,Tg,\dotsc,T^{r-1}g\}\sseq\cK$.
  By Lemma~\ref{lem-sublim}(a) and~(c),
  we either have $g=h$
  or $g,Tg,\dotsc,T^{r-1}g\lhd h$.
  But the latter option entails
  the contradictory inequality chain
  \[
    \stride g < \area g = \stride Tg < \dotsm 
    < \area T^{r-1}g = \stride T^rg = \stride g
  \]
  according to Corollary~\ref{cor-stride}.
\end{proof}

\begin{cor}\label{cor-unique}
  The IDE~\eqref{eqn-ide} has a solution only for $\gamma=\kappa$,
  and $h$ is the only unit stribola.
  Furthermore, the function $\ti h$ from~\eqref{eqn-tih}
  is the only standard stribola.
\end{cor}
\begin{proof}
  The assertions concerning~\eqref{eqn-ide} and~$h$
  follow immediately from Proposition~\ref{prp-fixed}
  and Theorem~\ref{thm-unique}.
  As for the standard stribola, suppose that $a>0$
  and $\ti g\colon[0,a]\to[0,a]$ satisfies $-\ti g'=\ti g^*$.
  Then $g:=\frac1a\ti g\mul a$ satisfies~\eqref{eqn-ide}
  with $\gamma=\frac1a$. Hence $g=h$, $a=\frac1\kappa$
  and $\ti g=\ti h$.
\end{proof}

We are now also ready to prove global convergence.

\begin{thm}\label{thm-global}
  We have $\lim_{n\to\infty}T^nf = h$ for every $f\in\cE$.
\end{thm}
\begin{proof}
  Let $f\in\cE$.
  In view of Lemma~\ref{lem-sublim}(d),
  it suffices to show that $h \in \cL(f)$.
  To this end, let us assume that $h\neq g_1\in\cL(f)$.
  Then, by Lemma~\ref{lem-sublim}(a)--(c)
  and Corollary~\ref{cor-stride},
  \begin{align}
    g_m &:= T^{m-1}g_1\in\cL(f), \label{eqn-gmLf}\\
    g_m &\lhd h, h_n \label{eqn-gmhhn}
    \quad\text{and}\\
    \gamma_{m-1} &:= \stride g_m 
     < \area g_m = \gamma_m = \stride g_{m+1}
     < \kappa \label{eqn-gmkappa}
  \end{align}
  for all $m,n\in\N$.
  In particular,
  $\lim_{m\to\infty}\frac{\gamma_m}{\gamma_{m-1}}=1$,
	hence, for each $n\in\N$ there exists $m\in\N$ such that
	$\frac{\gamma_m}{\gamma_{m-1}}\leq\frac{\kappa_n}\kappa$
	and thus
	$\kappa_{n+1}-1+\frac\kappa{\kappa_n} < \gamma_{m+1}$
	according to~\eqref{eqn-gmhhn}
	and Corollaries~\ref{cor-kappan} and~\ref{cor-hammock}.
	With \eqref{eqn-gmhhn} and~\eqref{eqn-gmkappa},
	we conclude
	that $(\gamma_m)_{m\in\N}$ converges to $\kappa$
	and that 	$(g_m)_{m\in\N}$ converges to $h$,
	implying that $h\in\cL(f)$ by~\eqref{eqn-gmLf}.
\end{proof}

The Hammock Lemma~\ref{lem-hammock}
also enables us to estimate $\kappa$ from below.

\begin{prp}
  For all $n\in\N$, we have
  $\kappa_n-1+\frac{\kappa_n}{\kappa_{n-1}} < \kappa$.
  In particular,
  \[
    0.2788770612338
    < \kappa_{23}-1+\frac{\kappa_{23}}{\kappa_{22}}
    < \kappa <  \kappa_{23} < 0.2788770613941.
  \]
\end{prp}
\begin{proof}
  Let $n\in\N$.
  The first assertion holds trivially for $n=1$
  and follows for $n>1$ by applying
  Hammock Lemma~\ref{lem-hammock}
  with $\ti f = h_n$, $f := h_{n-1}$ and $g := h$
  while using~\ref{cor-hdom}(a)--(b)
  and observing that
  $c := \frac{{\int}f}{{\int}g}
      = \frac{\kappa_{n-1}}{\kappa}
      = \frac{\stride h_n}{\stride h}$
  implies $h_n\leq h\mul\frac1c$ due to
  the Stride Lemma~\ref{lem-stride}.
  The values for $\thru\kappa0{23}$ are determined
  in~\cite{Mi}.
\end{proof}

\end{document}